\documentclass[12pt]{amsart}
\usepackage{amsfonts}
\usepackage{amsmath}
\usepackage{amsthm}

\newtheorem{thm}{Theorem}[section]
\newtheorem{prop}[thm]{Proposition}
\newtheorem{lem}[thm]{Lemma}
\newtheorem{hyp}[thm]{Hypothesis}
\newtheorem{rmk}[thm]{Remark}
\newtheorem{defn}[thm]{Definition}
\numberwithin{equation}{section}

\def\Or{\mathcal{O}_\mathbb{R}}
\def\Br{\mathbb{R}}
\def\Bc{\mathbb{C}}
\def\Co{\mathcal{O}}
\def\Xr{X_{\mathbb{R}}}
\def\Yr{Y_{\mathbb{R}}}
\def\c{\text{conj}}
\def\PN{P_N \langle \alpha \rangle}
\def\mfU{\mathfrak{U}}
\def\mfV{\mathfrak{V}}
\def\mfB{\mathfrak{B}}
\def\dbar{\overline{\partial}}
\def\Ca{\mathcal{A}}

\def\End{\text{End}}

\begin{document}
\title{A Dolbeault Isomorphism Theorem in Infinite Dimensions}
\author{Scott Simon}
\thanks{Research partially supported by NSF grant DMS 0203072.
I am especially grateful to Professor Lempert for his suggestions and patience.
This article contains the fruits of my thesis research under his guidance.
Also helpful were the suggestions of Aaron Zerhusen and
Parsa Bakhtary.}

\begin{abstract}
For a large class of separable Banach spaces, we prove the real analytic
 Dolbeault Isomorphism Theorem for open subsets.
\end{abstract}
\maketitle

\section{Introduction}\label{se:intro}

Dolbeault's isomorphism theorem states that if $E\rightarrow M$ is a
finite rank holomorphic vector bundle over a finite dimensional complex
manifold, then its sheaf and $\dbar$-cohomology groups are isomorphic:
\[
H^q(M,E) \approx H_{\dbar}^{0,q} (M,E).
\]
The case where $E$ is the trivial bundle can be found in \cite{Dolbeault}.
Our goal here is to extend this theorem to infinite dimensions.  An obvious
extension fails, even when $M$ is a domain in a Banach space.  Indeed, in
\cite{Patyi}, Patyi gives an example of a complex Banach space $X$
(which even has
an unconditional basis, see below) whose unit ball $B$ has
\[
H^q(B, \Co)=0, \qquad q \ge 1,
\]
but there is a closed $f \in C_{0,1}^\infty (B)$ that is not exact, hence
$H^{0,1} (B) \ne 0$.  We shall, however, show that a Dolbeault--type
isomorphism theorem can be proved in open sets in rather general Banach spaces
--- in particular in the Banach spaces Patyi considers --- if the Dolbeault
groups are defined in terms of real analytic forms.

Thus, let $X$ be a complex Banach space, $\Omega \subset X$ open,
$E \rightarrow \Omega$ a holomorphic Banach bundle, and
$\Ca_{p,q}=\Ca_{p,q}^{E}$ the sheaf of real analytic
$(p,q)$-forms on $\Omega$, with values in $E$, $p,q=0,1,2,\ldots$.
Then the operator
$\dbar : \Ca_{p,q} \rightarrow \Ca_{p,q+1}$
can be defined, much as in finite dimensions, see e.g.
\cite{Lempert:DolbeaultI}.
\begin{thm}\label{th:iso}
If $X$ has an unconditional basis, then
\begin{equation}\label{eq:iso}
H^q(\Omega,E) \approx
\frac{\text{Ker }\{\dbar : \Gamma (\Omega, \Ca_{0,q})   \rightarrow
\Gamma (\Omega, \Ca_{0,q+1})\}}
     {\text{Im } \{\dbar : \Gamma (\Omega, \Ca_{0,q-1}) \rightarrow
\Gamma (\Omega, \Ca_{0,q})  \}}.
\end{equation}
\end{thm}

Not surprisingly, the theorem will be obtained by considering the sheaf
$\Co^E$ of germs of holomorphic sections of $E$ and the complex
\begin{equation}\label{eq:complex}
0 \rightarrow \Co^E \rightarrow
\Ca_{0,0} \stackrel{\dbar}\rightarrow
\Ca_{0,1} \stackrel{\dbar}\rightarrow \ldots.
\end{equation}
It is known \cite[Proposition 3.2]{Lempert:DolbeaultI} that (\ref{eq:complex})
is exact, unlike its
$\mathcal{C}^\infty$ counterpart.  Therefore the abstract de Rham Theorem
(see, e.g., \cite[Theorem 3.13]{Wells})
would give (\ref{eq:iso}) if we knew that the sheaves $\Ca_{0,q}$
are acyclic, i.e. $H^p (\Omega, \Ca_{0,q})=0$ for $p \ge 1.$
This is what we are going to show, in fact in somewhat greater generality.
\begin{thm}\label{th:acyclic}
Let $\Xr$ be a real Banach space with unconditional basis,
$\Omega \subset \Xr$ open, $F \rightarrow \Omega$ a real analytic Banach
bundle, and $\Ca^r$ the sheaf of real analytic $r-$forms with values
in $F$.  Then $H^p (\Omega, \Ca^r)=0$ for $p \ge 1, \, r \ge 0.$
\end{thm}

In \cite{Cartan}, Cartan obtained a similar result in finite
dimensions.  As there, the key will be the following theorem:
\begin{thm}\label{th:nbhdbasis}
If $\Xr$ is a real Banach space with unconditional basis, $X \supset \Xr$
its complexification, then any set $S \subset \Xr$ has a neighborhood basis
in $X$ consisting of pseudoconvex open sets.
\end{thm}
The finite dimensional case follows a similar outline.  The key step was
a cohomology vanishing theorem which lead to the acyclicity of the resolution
\ref{eq:complex}.
In 1957, Cartan discussed 
the real analytic cohomology of real analytic manifolds \cite{Cartan}.
  If there is a
real analytic totally real imbedding into a complex manifold and the image
has a Stein neighborhood basis, then the corresponding
Theorems A and B for the sheaf of germs of real analytic sections
hold.  As noted above, the abstract de Rham Theorem, together with
acyclicity (Theorem B), imply the cohomology isomorphism theorem.
In 1958, Grauert completed the picture by proving the
necessary imbedding theorem as part of an
investigation of the Levi Problem \cite{Grauert}.  In Grauert's approach, a
pseudoconvex neighborhood basis again plays a
key role.  Grauert's results do 
not simply carry over to infinite dimensions,
since they rely on
compact sets with nonempty interior.  Such compact
sets are not available in infinite dimensions.  
Still, a part of the proof of Theorem~\ref{th:nbhdbasis}, namely the proof
of Theorem~\ref{th:nbhd} has some similarities 
to Grauert's method.

\section{Background}\label{se:bkgrnd}

Let $(X_\mathbb{R}, \|\,\|)$ be a real Banach space.
We define the complexification $X = \Xr \oplus i\Xr$ of $\Xr$ as this direct
sum of vector spaces, with the usual complex multiplication: if
$x_1,x_2 \in \Xr, \alpha, \beta \in \Br,$ then
\[
(\alpha + i\beta)(x_1 +ix_2)=
(\alpha x_1 - \beta x_2) + i(\beta x_1 + \alpha x_2).
\]
Given $x_1, x_2 \in \Xr$ and $x=x_1+ix_2,$ we define the projections
$\Re : X \rightarrow \Xr$ and
$\Im : X \rightarrow \Xr$ by $\Re x=x_1$ and
$\Im x=x_2.$
We define the norm
\[
\|x\|'=\sup_{0 \le \theta <2\pi}
\|\Re(e^{i\theta}x)\|.
\]
Since $\|\,\|$ agrees with $\|\,\|'$ on $\Xr,$ we will write
$\|\,\|$ for both.  Both $\Re$ and $\Im$ are real
linear maps of norm 1.  Define conjugation, a real linear isometry
 conj:$X\rightarrow X$, where $\c (x)=\Re x-i\Im x$.  We will
also write $x \mapsto \overline{x}$ for conjugation.  In general, when 
discussing the complexification of a Banach space, we will
use a Roman capital letter with the subscript $\Br$ to denote a real Banach 
space, and use the same letter
without the subscript to denote its complexification.  However, we will
sometimes denote a real Banach space by a Roman capital letter and use the
same letter with the subscript $\Bc$ to denote its complexification.
For more information on complexified Banach spaces, see \cite{Singer}.

Again, let $\Xr, \Yr$ be real Banach spaces.
Let $\Omega \subset X_\mathbb{R}$ be an open set,
 and let $f:\Omega \rightarrow Y_\mathbb{R}$.  Then $f$ is said to be real
 analytic if there are a neighborhood $U=\c (U) \subset X$ of $\Omega$ and
 a holomorphic function $g:U\rightarrow Y$ such that $g$ restricts to $f$ on 
$\Omega$. As in the finite dimensional case,
and in the differential and complex categories,
a real analytic Banach manifold is
a Hausdorff space glued together from open subsets of a Banach space
by real analytic ``gluing functions''.
A real analytic vector bundle is given by a real analytic map 
$\pi:E \rightarrow M$ of analytic Banach manifolds; each fiber of $\pi$ is
endowed with the structure of a real vector space.  It is required that for
each $x \in M $ there should exist a neighborhood $U \subset M,$ a Banach
space $\Yr$, and a real  analytic map $\pi^{-1} U \rightarrow U \times \Yr$
that has a real analytic inverse, and maps the fibers $\pi^{-1}$ linearly
on $\{\xi\} \times \Yr \approx \Yr, \xi \in U.$
A real analytic $E$-valued differential $r$-form
\[
f:\bigoplus^r TM \rightarrow E
\]
is a real analytic map which is multilinear on fibers, and respects 
foot-points.
This is all as in the finite dimensional case.
For precise definitions of manifolds, vector bundles, etc., see
\cite{Lempert:DolbeaultI}.

We wish to show that holomorphic functions are real analytic.  But first,
we will need some further background (found in \cite{Mujica}).
Given complex Banach spaces $W,Z$, an open set $\Omega \subset W$, and a
function
$f \in \Co (\Omega, Z),$ for each $a\in \Omega$ there are
$k$-homogeneous polynomials $P_k :W\rightarrow Z$ such that
\begin{equation}\label{eq:Taylor}
f(x)=\sum_{k=0}^\infty P_k(x-a)
\end{equation}
in a neighborhood of $a$. To each $k$-homogeneous
polynomial $P_k$ there is associated a unique symmetric $k$-linear map
$A_k$ such that $P_k(x)=A_k(x,\ldots x)$.
The sum \eqref{eq:Taylor} is a Taylor series, and the $k$th multiliniar map 
is proportional to the iterated $k$th differential of $f$.

Define
\[\|P_k\|=\sup_{\|x\|\le 1}\|P_k (x)\|
\]
and
\[
\|A_k\|=\sup_{\|x_1\|,\ldots \|x_k\| \le 1} \|A_k(x_1,\ldots x_k)\|.
\]
Then $\|A_k\|\le e^k\|P_k\|$ by \cite[Exercise 2.G]{Mujica}.
The Cauchy-Hadamard formula tells us
that the radius of uniform convergence, $R$, of (\ref{eq:Taylor}) is given by
$R^{-1}=\limsup_k \|P_k\|^{1/k}.$
\begin{prop}\label{uniqueness}
Let $X$ and $Y$ be complexified Banach spaces, $G\subset X$ open and connected,
$G\cap \Xr \ne \emptyset,$ and $f\in \Co (G,Y).$  If $f|_{G\cap \Xr}=0,$ then
$f=0$.
\end{prop}
\proof
Choose $a\in G\cap \Xr,$ and write
\[
f(x)=f(a)+\sum_{k=1}^\infty A_k (x-a,\ldots x-a).
\]
By assumption, $f(a)=0.$  But since $A_k$ is proportional to the $k$th
differential of $f,$ $A_k=0$ in all real directions.  Since any $x\in X$
is a linear combination of $\Re x +i \Im x,$ $A_k=0$ by complex multilinearity.
Therefore, $f=0.$
\endproof

Now we can prove that holomorphic functions are real analytic.
\begin{prop}
Given complex Banach spaces $W,Z$, and an open set $\Omega \subset W$,
then $\Co (\Omega, Z) \subset \Ca (\Omega,Z).$
\end{prop}
\proof
For $a\in \Omega$ arbitrary, write
\begin{equation}\label{eq:Taylor2}
f(x)=\sum_{k=0}^\infty P_k(x-a)=\sum_{k=0}^\infty A_k (x-a,\ldots x-a).
\end{equation}
Since any complex
multilinear map is also real linear, we can disregard the complex structure on
$W$ and $Z$, and still write $f$ as a sum of real multilinear maps
near $a$.
We may extend $A_k$ to the complexification $W_{\Bc}$ of $W$ in the natural
way; if $x_1,\ldots x_k,y_1, \ldots y_k \in W,$ and multiplication by $i$
refers to the complex multiplication in $W_\Bc$ and $Z_\Bc$, define the 
extension
$\tilde{A}_k:W_{\Bc}\times \ldots \times W_{\Bc} \rightarrow Z_{\Bc}$
 of $A_k$ by the
formula
\begin{multline*}
\tilde{A}_k(x_1+iy_1, \ldots, x_k+iy_k)=A_k(x_1,\ldots, x_k)+
iA_k(x_1,\ldots ,x_{k-1}, y_k)\\
+iA_k(x_1, \ldots, x_{k-2},y_{k-1}, x_k)
+\ldots +iA_k(y_1,x_2,\ldots, x_k)\\
+i^2A_k(x_1,\ldots, x_{k-2},y_{k-1},y_k)
+\ldots +i^2A_k(y_1,y_2,x_3,\ldots, x_k)+\ldots \\
+i^kA_k(y_1,\ldots, y_k).
\end{multline*}
Observe that there are $2^k$ terms in this sum.  To each $\tilde{A}_k$,
a homogeneous polynomial $\tilde{P}_k$ can be associated.  We now investigate
the convergence of
\begin{equation}\label{eq:Taylor3}
\sum_{k=0}^\infty \tilde{P}_k ((x+iy) - a).
\end{equation}
We can estimate the terms:
\begin{equation*}
\begin{split}
\|\tilde{P}_k ((x+iy) - a)\| 
&\le 2^k\|A_k\|\max(\|x-a\|,\|y\|)\\
&\le 2^ke^k\|P_k\|\max(\|x-a\|,\|y\|)^k.
\end{split}
\end{equation*}
According to the Cauchy-Hadamard formula, the radius of uniform convergence
$R$ of \eqref{eq:Taylor2} is given by
\[
R^{-1}= \limsup_{k} \|P_k\|^{1/k}.
\]
Thus, if we use the root test for the sum
\[
\sum_{k=0}^\infty \|\tilde{P}_k ((x+iy) - a)\|,
\]
we have
\begin{align*}
\limsup_k \|\tilde{P}_k ((x+iy) - a)\|^{1/k}\\
\le 2e\max(\|x-a\|,\|y\|)\limsup_k \|P_k\|^{1/k}\\
\le 2eR^{-1}\max(\|x-a\|,\|y\|).
\end{align*}
Therefore \eqref{eq:Taylor3} converges uniformly on a neighborhood containing 
the set
\[
\{x+iy \in X:\max(\|x-a\|,\|y\|)<R/2e\}.
\]
Then $f$ can be extended to a holomorphic
function in a $W_{\Bc}$-neighborhood of any $a\in \Omega.$  But we would
like to find a single neighborhood $G \subset W_{\Bc}$ of $\Omega$ to which
we can extend $f.$  We will check that the sums
\eqref{eq:Taylor3} which converge in a ball centered at each $a \in \Omega$
agree on the overlaps of these balls.  If so, then they
define a holomorphic function on a neighborhood of $\Omega.$
Suppose $g$ and $g'$ are
two extensions of $f$ on two overlapping
balls.  Then both agree on the intersection of their domains
intersected with $\Omega.$  By Proposition~\ref{uniqueness}, $g-g'=0$ on
the entire intersection of their domains, or $g=g'$ on the overlap, as
required.
\endproof
\begin{defn}
Given complexified Banach spaces $X$ and $Y$, and an open set $U\subset X$
such that $U=\c (U)$, a
 function $f:U\rightarrow Y$ is called ``real-type holomorphic'' or, ``of real
type'' (written $f\in \Or (U,Y)$) if $f$ is holomorphic and commutes with
conjugation, i.e. $\overline{f(x)} = f(\overline x)$.
\end{defn}
If this is the case, and if $x \in U \cap \Xr$, then
$f(x) \in Y_{\mathbb{R}}.$ Observe that the sum,
direct sum, composition, etc. of two real-type holomorphic
functions is again of real type.  If a sequence of real-type holomorphic
functions converges locally uniformly, the limit is also a real-type
holomorphic function.

A real (complex) Banach space $X$ has a  Schauder basis if there are
$\{e_j\}_{j\in \mathbb{N}} \subset X$ such that any $x \in X$ can be
written uniquely as a sum
\begin{equation}\label{eq:Schauder}
\sum_{j=1}^\infty \lambda_j e_j, \qquad \lambda_j \in \Br\, (\Bc).
\end{equation}
If $(\ref{eq:Schauder})$ converges unconditionally, i.e. independently of
any rearrangement of terms for all $x \in X,$ then
$\{e_j\}$ is said to be an unconditional basis.  A Schauder basis
of a real Banach space $\Xr$ is also a Schauder basis of its complexifictation,
$X,$ and an unconditional basis
of $\Xr$ is also an unconditional basis of $X$.

If $X$ is a complex Banach space, a pseudoconvex set $P\subset X$ is an open 
set such that $-\log (\text{dist}(x,\partial P))$ is plurisubharmonic in $P$.
Observe that $X$ itself is also pseudoconvex.

\section{A Pseudoconvex Neighborhood Basis}
Theorem~\ref{th:nbhdbasis} is proved in two steps.  The first step is the
following theorem.

\begin{thm}\label{th:nbhd}
Let $X$ be a complexified separable Banach space and let $\Omega \subset
 X_\mathbb{R}$ be open.  Then there is a pseudoconvex $P \subset  X$ such that
 $P \cap X_\mathbb{R} = \Omega$.
\end{thm}
\proof
If $\Omega = \Xr,$ then $P = X$ will suffice.  Otherwise,
$\Xr \setminus \Omega$
is nonempty.
Let $A \subset \Omega$  be a countable dense subset.
 For every $a\in A,b\in \Xr \setminus \Omega,$ define
 $l_{ab}$ to be a real linear
functional on $\Xr$ of norm 1 such that $l_{ab} (a-b)=\|a-b\|.$
Extend $l_{ab}$ to a complex linear functional of real type by setting
$l_{ab}x=l_{ab}(\Re x) +il_{ab}(\Im x)$.
Now let
$\{a_k\}_{k=0}^\infty$ be an enumeration of $A$, and set
\begin{align}\label{eq:phib}
\phi_b (x)&=\sum_{k=0}^\infty  4^{-k}(l_{a_k b} (x-b))^2 \\
\notag
&=\sum_{k=0}^\infty  4^{-k}\left((l_{a_k b} x)^2-2(l_{a_k b} x)(l_{a_k b} b)
+(l_{a_k b} b)^2\right),
\end{align}
whence for any $M>0,$ the family $\{\phi_b\}, \|b\| \le M$ is uniformly
equicontinuous on bounded subsets of $X$.
Since the series is locally uniformly convergent,
$\phi_b$ is an entire holomorphic function.  As the sum of real--type 
functions, $\phi_b$ is itself of real type.
Furthermore, on $\Omega$ we can estimate
$\phi_b$ from below uniformly in $b$:  Let
$r(x)=\text{dist}(x,\Xr \setminus \Omega),$  and fix $x \in \Omega.$
Choose $a_{k_0}\in A$ such that $\|x-a_{k_0}\|<r(x)/4.$  Observe that
\begin{align*}
 \sum_{k=0}^\infty 4^{-k}(l_{a_k b} (x-b))^2
 &\ge 4^{-k_0}(l_{a_{k_0} b} (x-b))^2\\
 &\ge 4^{-k_0}(l_{a_{k_0} b} (a_{k_0} -b)-l_{a_k b} (x-a_{k_0}))^2\\
 &\ge 4^{-k_0}((\|a_k-b\| - \|x-a_{k_0}\|)^2\\
 & >  4^{-k_0}(r(a_{k_0})-r(x)/4)^2 \\
 &\ge 4^{-k_0}(3r(x)/4-r(x)/4)^2\\
 &\ge 4^{-k_0-1}r(x)^2>0.
\end{align*}
Thus, if we set
\[
u(x) =\sup_{b \in \Xr \setminus \Omega} -\Re (\phi_b (x)), \qquad x\in X
\]
we have $u(x) < 0$ for $x \in \Omega.$ Furthermore,
$u(x)=0$ on $\Xr \setminus \Omega,$ since clearly $u \le 0$ on $\Xr,$
and for any $b \in \Xr \setminus \Omega,$ $\phi_b (b)=0.$

Next we show that $u$ is continuous.  From \eqref{eq:phib},
\begin{align*}
&\liminf_{\|b\|\rightarrow \infty} \Re \phi_b(x)/\|b\|^2\\
&\ge \liminf_{\|b\|\rightarrow \infty}
\left\{\Re (l_{a_0b}(a_0-b)+l_{a_0b}(x-a_0))^2
-\sum_{k=1}^\infty 4^{-k}\|x-b\|^2\right\}/\|b\|^2\\
&\ge \liminf_{\|b\|\rightarrow \infty} \|a_0-b\|^2/\|b\|^2-1/2 \ge 1/2,
\end{align*}
uniformly for $x$ in a bounded set $V$.  Hence, given $V,$ 
$\Re \phi_b \rightarrow \infty$ as $\|b\| \rightarrow \infty$ uniformly on $V,$
and so 
\[u(x)=\sup_{b \in \Xr \setminus \Omega} -\Re \phi_b(x)
=\sup_{b \in \Xr \setminus \Omega, \|b\| \le M} -\Re \phi_b(x)
\]
for all $x \in V$, provided $M$ is sufficiently large.  Since
$\{\phi_b:\|b\| \le M\}$ is equicontinuous on bounded subsets, it follows that
$u$ is continuous.  As a continuous supremum of 
plurisubharmonic functions,
it is also plurisubharmonic.   Therefore $P=\{x\in X:u <0\}$ is a
pseudoconvex neighborhood of $\Omega.$  Since $u=0$ on
$\Xr \setminus \Omega$ and $u<0$ on $\Omega,$ $P \cap \Xr = \Omega.$
\endproof
The second step in the proof of Theorem~\ref{th:nbhdbasis} is the heart of the
 entire matter. The goal is to show that arbitrarily ``narrow" pseudoconvex
 neighborhoods exist. The critical tool is a theorem about
 real-type holomorphic domination.  But
first, given a Banach space $X,$ $a \in X,$ $r \in \Br,$ define 
\[B(a;r)=\{x\in X: \|x-a\|<r\}\] and
$B(r)=B(0;r).$  We state a Runge-type hypothesis which we will use in
the following theorem:
\begin{hyp}\label{hyp:Runge}
There is a $\mu\in (0,1)$ such that for any Banach
space $(W,\| \ \|_{W})$, $\epsilon >0$, and $g\in \Co (B(1);W)$
there is an $h\in \Co (X;W)$ that satisfies $\| g-h\|_W<\epsilon$ on
$B(\mu)$.
\end{hyp}
This allows us to state the critical theorem.
\begin{thm}\label{th:RHD}
Let $X$ be a complexified Banach space with a Schauder Basis
satisfying Hypothesis~\ref{hyp:Runge}, let
$P \subset X$ be pseudoconvex, and suppose $P=\c (P)$.  Let
$u :P\rightarrow \mathbb{R}$ be locally Lipschitz, with $u \le 0$ on
$P\cap X_\mathbb{R}$.
 Then there is a complexified Banach space $Y$ and an $f\in \Or (P,Y)$
such that $u \le \|\Im f\|$.
\end{thm}
Lempert proved an analogous theorem in \cite{Lempert:PSHD}.
Using Theorem~\ref{th:nbhd} and Theorem~\ref{th:RHD} together, we can prove
the following generalized version of Theorem~\ref{th:nbhdbasis}:

\begin{thm}\label{th:gen_nbhdbasis}
If $\Xr$ is a real Banach space with a Schauder basis
and its complexification, $X \supset \Xr$, satisfies
Hypothesis~\ref{hyp:Runge}, then any set $S \subset \Xr$ has a
neighborhood basis in $X$ consisting of pseudoconvex open sets.
\end{thm}

This theorem implies Theorem~\ref{th:nbhdbasis} because a Banach space
with an unconditional basis satisfies Hypothesis~\ref{hyp:Runge} by
\cite{Lempert:ApproxII}.

\proof[Proof of Theorem~\ref{th:gen_nbhdbasis}]
Let $G \subset X$ be an arbitrary open neighborhood of $S.$  We 
construct a pseudoconvex neighborhood of $S$ contained in $G.$  By
Theorem~\ref{th:nbhd}, there is a pseudoconvex neighborhood $P$ of $G \cap \Xr$
such that $P \cap \Xr = G \cap \Xr.$  By passing to $P \cap \c (P),$ we can 
asume
$P = \c (P).$  We can also assume $G \subset P.$
Define $u: P \rightarrow \Br$ by
\[
 u (x)=
 \left\{
  \begin{array}{ll}
   \min (1,\|\Im x\|/\text{dist}(x,\partial G)),
    & \text{if } x\in G,\\
   1& \text{otherwise.}
  \end{array}
 \right.
\]
$u$ is locally Lipschitz: in $P\setminus \overline{G}$ this is obvious, and
in $G$ it follows from the fact that
$\text{dist}(x,\partial G)$ is Lipschitz.  Now 
$(P \setminus G) \cap \Xr= \emptyset,$ so if $x \in P \cap \partial G$
then  $\|\Im  x\| >0.$  Therefore dist$(y,\partial G) < \|\Im y\|$ for
$y$ in some neighborhood of $x$, and so $u=1$ is Lipschitz there too.
Theorem~\ref{th:RHD} implies that there is a complexified Banach space
$Y$ and an $f\in \Or (P,Y)$ such that $u \le \|\Im f\|$.
Set $Q=\{\|\Im f\| < 1\}.$  Then $Q$ is pseudoconvex, and  
$S \subset P \cap \Xr \subset Q$ since $\Im f = 0$
on $\Xr$.  Furthermore, if $x\in Q,$ then
$\|\Im f(x)\| <1,$ so $u(x) < 1,$ which implies that
$x \in G$, as required.
\endproof

Theorem~\ref{th:RHD} remains to be proven.  As in the proof of an analogous
theorem in \cite{Lempert:PSHD},
this proof will be by induction.
Given an open set $P\subset X,$ consider those balls
$B=B(a;r)$ such that
\renewcommand\theenumi{\roman{enumi}}
\renewcommand\labelenumi{(\theenumi)}
\begin{enumerate}
  \item $\overline{B(a;r)} \cup \c (\overline{B(a;r)})\subset P,$
  \item $2\text{diam } B \le \text{diam } P,$
  \item $B= \c (B)\text{ or }\overline{B} \cap \overline{\c (B)}
    =\emptyset.$
\end{enumerate}
Let $\mfB_P$ denote this family of balls.

Now we are ready to formulate the induction step in the form of the
following proposition.
\begin{prop}\label{pr:descent}
Let $X$ be a complexified Banach space with Schauder basis satsifying
Hypothesis~\ref{hyp:Runge}.  Let $P=\c (P) \subset X$ be a
pseudoconvex set.
If for every $B\in \mfB_P$
there are a complexified
 Banach space $(V_B, \|\,\|_B)$ and an $f_B\in \Or (B\cup \c (B),V_B)$
 such that
 $u \le \|\Im f_B\|_B$ on $B$, then there is a complexified Banach space 
$(V,\|\,\|_V)$
and an $f\in \Or (P,V) $ such that
 $u \le \|\Im f\|_V$ on $P$.
\end{prop}
This proposition implies Theorem~\ref{th:RHD}.
Its proof will take up sections~\ref{se:ballbndl} and~\ref{se:descent}.

\proof[Proof of Theorem~\ref{th:RHD}]
Suppose not.  If $u$ cannot be dominated on $P$ by any $\|\Im f \|$ such
that $f\in \Or (P,Y)$,  then
by Proposition~\ref{pr:descent}, there is a $B_1 \in \mfB_P$ such that 
$u$ cannot be dominated in the same way on
$B_1 \cup \c(B_1)$.  Replacing $P$ with $B_1 \cup \c (B_1)$, we can repeat
the same
argument to produce $B_{k+1} \in \mfB_{B_k}$ such that 
$\overline{B_{k+1}} \subset B_k$
and $u$ cannot be dominated on $B_{k+1} \cup \c (B_{k+1}).$  Since 
2diam $B_{k+1} \le \text{diam }B_k,$ the $B_k$ converge to a point $x_0 \in P.$
Choose $r>0$ such that 
$B=B(x_0;r) \in \mfB_P,$ and $u$ has some Lipschitz constant $K>0$ on 
$B \cup \c (B)$. Since $B_k \subset B$ for some $k,$ $u$ cannot be dominated
on $B \cup \c (B)$.  

Suppose first that $\overline{B} \cap \Xr = \emptyset$.
Then $\Re x_0 \notin \overline{B},$ so that \linebreak
$\|x_0-\Re x_0\|>r.$  Set
\[
f(x)=\frac{|u(x_0)|+|u(\overline{x_0})|+Kr}{\|\Im x_0\|-r} x.
\]
Since
$\|\Im x\| \ge \|\Im x_0\|-r$ when $x \in B,$ we have
$\|\Im f(x) \| \ge |u(x_0)|+Kr \ge u(x)$; and similarly $\|\Im f \| \ge u$
on $\c(B).$

In the second case, $B=\c (B)$. With $f(x)=Kx,$
\[
u (x) \le |u (\Re x)| + K\|x-\Re x\| \le 0+K\|\Im x\|=\|\Im f(x)\|.
\]
Thus in either case there is an $f \in \Co (B \cup \c(B),\Bc)$ such that
$\|\Im f\| \ge u,$ contradicting our earlier finding.  This contradiction
proves the theorem.
\endproof

\section{Ball Bundles}\label{se:ballbndl}

One of the tools which we will use to prove Proposition~\ref{pr:descent}
is ball bundles over finite dimensional bases.  The setup
is exactly the same as in, e.g. \cite{Lempert:PSHD}.  For convenience, we 
include the necessary definitions and propositions here. 
After renorming as in \cite[Section 7]{Lempert:DolbeaultIII}, we may assume
that with a given basis $\{e_j\}$ of $\Xr$
and $1 \le n \le N \le M \le m \le \infty,$
\[
\|\sum_{j=N}^M \lambda_j e_j \| \le \|\sum_{j=n}^m \lambda_j e_j \|, \qquad
\lambda_j \in \Bc.
\]
This renorming respects conjugation, etc.

Let $\pi_N$ be a projection on the first $N$ coordinates, and 
$\rho_N=id-\pi_N.$
Fix $P \subset X$ pseudoconvex.
Let  $d(x) = \min\{1, \text{dist}(x,X\setminus P) \}$ and, given
$0 < \alpha < 1,$
\begin{align*}
D_N \langle \alpha \rangle &= \{ t \in \pi_N X : \|t \| < \alpha N, 1
< \alpha N d(t) \}, \\
P_N \langle \alpha \rangle &= \{x \in X : \pi_N x \in D_N \langle \alpha
 \rangle, \| \rho_N x
 \| < \alpha
 d(\pi_N x) \}.
\end{align*}
 These sets have the following properties (proved in \cite{Lempert:PSHD}):
\begin{prop}\label{pr:ballprops}
For any pseudoconvex set $P$, integer $N,$ and number $\alpha$ with 
$0<\alpha<1,$ the following hold:
\renewcommand\theenumi{\alph{enumi}}
\renewcommand\labelenumi{(\theenumi)}
\begin{enumerate}
 \item Each $P_N \langle \alpha \rangle \subset P$ is pseudoconvex.
 \item For fixed $\alpha,$ each $x \in P$ has a neighborhood that is contained
       in all but finitely many  
       $P_N \langle \alpha \rangle$.
\end{enumerate}
\end{prop}
We will use an approximation theorem for ball bundles:
\begin{thm}\label{th:ballapprox}
Assume Hypothesis~\ref{hyp:Runge} with some $\mu\in(0,1)$.  If
$\gamma<2^{-6}\mu\alpha$ and $V$ is a complex Banach space then any
$\psi\in \Co (P_N\langle\alpha\rangle;V)$ can be approximated by 
$\phi\in\Co (P;V)$, uniformly on $P_N\langle\gamma\rangle$.
\end{thm}
This theorem is the same as \cite[Theorem 3.3]{Lempert:DolbeaultIII}.

Here is a proposition relating ball bundles to $\mathfrak{B}_P$.

\begin{prop}\label{pr:cover}
Let $X$ be a complexified Banach space, and $P=\c (P) \subset X$ pseudoconvex.
For any positive integer $N$, and any choice of $\alpha$ satisfying
$0<2^7\alpha < \mu^2 <1,$ $\PN$ has
a finite cover by balls $B_k=B(x_k,r_k)$ such that
$B(x_k, 2r_k/\mu) \in \mfB_P.$
\end{prop}

\proof
Let $A=\overline{P_N\langle\alpha\rangle}\cap \pi_N X.$  Then
$P_N\langle\alpha\rangle\cap \pi_N^{-1}t \subset B(t;\alpha d(t))$ for any
$t\in A$.  For each $t$ there is a relatively open $U_t \subset \pi_N P$
such that
$P_N\langle\alpha\rangle\cap \pi_N^{-1}U_t \subset B(t;2\alpha d(t)).$
Since $\{U_t\}_{t\in A}$ covers $A$, we can find a finite set $T\subset A$
such that $\{U_t\}_{t\in T}$ covers $A$.  But then
$\{B(t;2\alpha d(t))\}_{t\in T}$ covers $\PN$.

We claim that each of these balls is contained in an element of 
$\mathfrak{B}_P.$
To prove the claim, fix $t \in A.$

Case 1: $4\alpha d(t)/\mu < \text{dist} (t,\Xr)$.\\
We check that
$B(t;4\alpha d(t)/\mu) \in \mathfrak{B}_P.$
Since
$4\alpha d(t)/\mu < d(t),$ we have
$ \overline{B(t;4\alpha d(t)/\mu) } \subset P.$
Also,
\[
 2\text{diam}(B(t;4\alpha d(t)/\mu))=16\alpha d(t)/\mu < 2d(t)
\le \text{diam} (P).
\]
Furthermore,
\[
B(t;4\alpha d(t)/\mu)\cap \c (B(t;4\alpha d(t)/\mu)) = \emptyset.
\]
Therefore,
\[
B(t;4\alpha d(t)/\mu) \in \mathfrak{B}_P.
\]

Case 2:
$4\alpha d(t)/\mu \ge \text{dist} (t,\Xr)$.
We will find a ball $B(s;\mu d(s)/4)$ with $s \in \Xr$ which contains
$B(t;2\alpha d(t))$, and then show that
\[
B(s;d(s)/2) \in \mathfrak{B}_P.
\]
Choose $s \in P \cap \Xr$ such that $\|s-t\| < 8\alpha d(t)/\mu.$
Let
$x\in B(t;4\alpha d(t)/\mu).$  Then
\[
\|x-s\| \le \|x-t\|+\|t-s\| < 16\alpha d(t)/\mu.
\]

But
\[
d(t) \le d(s) + \|s-t\| \le d(s) +8\alpha d(t)/\mu,
\]
so
\[
d(t) \le d(s)(1-4\alpha/\mu) <2d(s).
\]
Therefore, $\|x-s\| < 2^5\alpha d(s)/\mu<\mu d(s)/4.$  In other words,
\[
B(t;4\alpha d(t)/\mu) \subset B(s;\mu d(s)/4).
\]
We check that $B(s;d(s)/2) \in \mathfrak{B}_P.$
Since $s \in \Xr,$
$B(s;d(s)/2)= \c (B(s;d(s)/2).$
Clearly,
$\overline{B(s;d(s)/2)} \subset P.$
Furthermore,
\[2\text{diam}(B(s;d(s)/2) =2d(s) \le \text{diam}(P).\]
Then  $B(s;d(s)/2) \in \mathfrak{B}_P$.
\endproof

\section{The proof of Proposition~\ref{pr:descent}}\label{se:descent}

In this section, we will use the following conventions when taking direct
sums of Banach spaces.
Given Banach spaces $(X_j,\|\, \|_j), j\in J,$ define the Banach space
$\overline{\bigoplus}_{j\in J} X_j$ to be the set of all bounded collections
$(x_j), x_j\in X_j$ with the supnorm
$\|x\|=\sup \|x_j\|_j.$
Observe that the complexification of
$\overline{\bigoplus}_{j\in J} X_{j\Br}$ is
$\overline{\bigoplus}_{j\in J} X_j .$
\begin{rmk}\label{re:decomp}
Given complexified Banach spaces $X,Y$ an open $G\subset X$ such that
$G=\c (G)$ and a function $f\in \Co (G,Y),$ define
$f' (x)=(f(x)+\overline{f(\overline{x})})/2,$ and
$f''(x)=(f(x)-\overline{f(\overline{x})})/2i.$
Then $f',f''\in \Or (G,Y),$
$f'(x)\oplus f''(x)\in \Or (G,Y\oplus Y),$
 and 
\[\|f(x)\|\le 2\|f'(x)\oplus f''(x)\|\le 2\max(\|f(x)\|,\|f(\overline{x})\|).\]
Furthermore, if $f$ is a bounded linear map, then $\|f'\|,\|f''\| \le \|f\|,$
and
\[\|f(x)\| \le 2\max(\|f'(x)\|,\|f''(x)\|)\]
for all $x \in G.$
\end{rmk}

Some further results are required to prove Proposition~\ref{pr:descent}.
If $A\subset \pi_NX\approx \mathbb{C}^N$ and $r:A\to [0,\infty)$ is continuous,
 define the sets
\begin{align*}
A(r) &= \{x\in X: \pi_N x\in A, \|\rho_N x\| < r(\pi_N x)\}\\
A[r] &= \{x\in X: \pi_N x\in A, \|\rho_N x\|\le r(\pi_N x)\}.
\end{align*}
\begin{lem}\label{lem:tech1}
Let $X$ be a complexified Banach space with Schauder basis, $P =\c (P)
 \subset X$ pseudoconvex, $N$ a positive integer, $A_1 \subset \subset A_2
 \subset \subset A_3 \subset \subset A_4 \subset \pi_N P,$ with
$A_i=\c (A_i)$ open, $i=2,3,4,$ $A_1= \c (A_1)$
 compact and plurisubharmonically convex in $A_4$.
   Let $r_i:A_4 \rightarrow (0,\infty)$ be continuous,
    with $2r_1<r_2<r_3<r_4$,
 $-\log r_1$ plurisubharmonic on $A_4$, and
 $r_i (x)=r_i(\overline{x})$ for $x \in A_4.$   Finally, suppose that all 
Banach space valued 
holomorphic functions on $A_4(r_4)$
can be approximated by holomorphic
 Banach space valued functions on $P$, uniformly on $A_3(r_3).$  Then for any 
complexified Banach
 space $V$ and any function $f\in \Or (X,V),$ there is a complexified 
Banach space $W$ and a $g \in \Or (P,W)$ such that
\renewcommand\theenumi{\roman{enumi}}
\renewcommand\labelenumi{(\theenumi)}
\begin{enumerate}
 \item $\|g\| \le 1$ on $A_1[r_1]$, and
 \item $\|\Im g \| \ge \| \Im f \|$ on $A_3(r_3) \setminus A_2(r_2)$.
\end{enumerate}
\end{lem}
\proof
Set $\pi = \pi_N,$ $\rho = \rho_N$.
First, we produce a complexified Banach space $Z$ and a function
$\phi \in \Or (P,Z)$
 which has norm less than 1/8 on $A_1 [r_1]$
and greater than 2 on
$A_3(r_3)\setminus A_2(r_2).$
Consider the constant function equal to 4 on $X$.  By 
\cite[Lemma 4.1]{Lempert:PSHD}, there
is a complex vector space $V$ and a function $\phi_1 \in \Co (X,V)$ which has
norm less than 1/8 on $A_1 [r_1]$ and greater than 4 on
$A_3(r_3)\setminus A_2(r_2).$  Observe that the vector space $V$
produced in \cite[Lemma 4.1]{Lempert:PSHD} is
complexified.  Then by Remark~\ref{re:decomp}, we see that
$Z=V\oplus V,$ $\phi(x)=\phi_1'(x)\oplus \phi_1'' (x)$ will do.

Now we will need an auxiliary family of functions.
Define
\[
w_\lambda (z)=(\lambda^2- (z-\lambda)^2)= 2\lambda z-z^2, \qquad -1 \le 
\lambda \le 1.
\]
Then $w_\lambda \in \Or (\Bc)$.  Furthermore, $w_\lambda$ satisfies:
\renewcommand\theenumi{\roman{enumi}}
\renewcommand\labelenumi{(\theenumi)}
\begin{enumerate}
\item
$|w_\lambda (z)| \le 1/2$ whenever $|z| \le 1/8,\,1 \le \lambda \le 1,$
\item
$w_\lambda (z)=|z|^2$ whenever $\lambda = \Re z$, and
\item
$w_\lambda$ are uniformly bounded on bounded subsets of $\Bc$.
\end{enumerate}
We are ready to define $W$ and $g\in \Or(P,W)$.  First, let $K$
be the set of all real-type linear functionals in the closed unit ball
 of the dual $Z'.$  Let

 \[W=\{\text{bounded maps from }K \times [-1,1] \text{ to }\Bc \},\]
with the sup norm; it can be identified with the complexification of
 \[W_\Br=\{\text{bounded maps from }K \times [-1,1] \text{ to }\Br \}.\]
Choose
$q$  large enough so that
\[
\|f(x)\| \le 2^q \qquad \text{whenever }\pi x \in A_1 \text{ and}
\|\rho x\| \le 2^{-q} \max_{A_1} r_1.
\]
Let
\[
g(x)(k,\lambda)=
(w_\lambda (k\phi (x)))^q f(\pi x + (w_\lambda (k\phi (x)))^q \rho x).
\]
If $x \in A_1[r_1],$ then
$(w_\lambda (k\phi (x)))^q \le 2^{-q},$
 and
$\|f(\pi x + (w_\lambda (k\phi (x)))^q \rho x)\| \le 2^{q},$
so
$\|g(x)\| \le 1$.  On the other hand, if
$x \in A_3(r_3) \setminus A_2(r_2),$ then $\|\phi (x)\| \ge 2$, so the 
Hahn-Banach Theorem and Remark~\ref{re:decomp} imply that there is
a $k \in K$ such that $|k \phi (x)| = 1$.  If $\lambda = \Re (k \phi (x)),$
then $w_\lambda (k\phi (x))=1.$  In this case, $g(x)(k,\lambda)=f(x),$ so
$\|\Im g(x) \|\ge \|\Im f(x) \|$.  Furthermore, since $w_\lambda,k,\phi, \pi,$
and $\rho$ are of real type, so is $g.$
\endproof

Lemma~\ref{lem:tech1} implies the following proposition, similar to
\cite[Proposition 4.2]{Lempert:PSHD}:
\begin{prop}\label{pr:tech2}
Assume Hypothesis~\ref{hyp:Runge}, and let $2^5\beta<\alpha<2^{-8}\mu$.
If $Z$ is a complexified
 Banach space and $g\in \Or(X;Z)$, then there are a complexified Banach
 space $W$ and $h\in \Or(P;W)$ such that
\renewcommand\theenumi{\roman{enumi}}
\renewcommand\labelenumi{(\theenumi)}
\begin{enumerate}
\item $\| h(x)\|_W \le 1$ if $x\in P_N\langle\beta\rangle$, and
\item $\|\Im h(x)\|_W\ge\|\Im g(x)\|_Z$ if
$x\in P_{N+1}\langle\alpha\rangle\setminus P_N\langle\alpha\rangle$.
\end{enumerate}
\end{prop}

\proof
Recall that $d(x)=\min(\text{dist}(x,\partial P),1)$.
Note that since $2^5\beta<\alpha,$ if
$r_1= 4\beta d, r_2 =\alpha d/4,$ then $2r_1 < r_2$.  The rest of the proof
is the same as in \cite[Proposition 4.2]{Lempert:PSHD}, except substituting
Lemma~\ref{lem:tech1} for \cite[Lemma 4.1]{Lempert:PSHD}.
\endproof

\begin{prop}\label{pr:balldom}
Let $X, V$ be complexified Banach spaces such that Hypothesis~\ref{hyp:Runge}
 holds for $X$, and let $B = B(a;r) \subset X$ satisfy
$\overline{B} \cap \overline{\c (B)} = \overline{B}$ or $\emptyset.$
If $ U=B \cup \c (B)$ and
$f \in \Or (U,V),$  then there is a complexified Banach space $W$ and a 
function $g \in \Or (X,W)$ such that $\|\Im f \| \le \|\Im g \|$ on
$B(a;\mu r/2).$
\end{prop}

\proof
\item{Case 1: $B = \c (B)$.}  Then $a \in \Xr.$
By Hypothesis~\ref{hyp:Runge}, there is an $h \in \Co (X,V)$ such that
$\|f-h\| \le 1$ on $B(a;\mu r)$.  In fact, replacing $h$ by $h'$ as defined
in Remark~\ref{re:decomp}, we can assume $h \in \Or(X,V).$
It follows that $f-h$ and $\Im (f-h)$
are Lipschitz on the ball $B(a; \mu r/2)$, with some
Lipschitz constant $M$.  Let $W=V \oplus X$ with
the supnorm, and $g(x)= h(x) \oplus Mx.$  Then for any $x \in B(a;\mu r/2),$
\begin{align*}
\| \Im f(x) \| &\le \|\Im h(x)\| + \|\Im (f(x) - h(x))\| \\
               &\le \|\Im h(x)\| + \|\Im\{(f(x) - h(x))
                -(f(\Re (x))-h(\Re (x)))\} \| \\
               &\le \|\Im h(x)\|+M\|\Im x\|=\|\Im g(x)\|.
\end{align*}

\item{Case 2: $\overline{B} \cap \overline{ \c (B)} = \emptyset$} (and
therefore $\|\Im a \| -r >0$).
By Hypothesis~\ref{hyp:Runge}, there is a function
$g \in \Co (X,V\oplus \Bc)$ such that
\[
\|g - (2f \oplus 2x/(\|\Im a \| -r)) \| <1 \qquad \text{on } B(a; \mu r).
\]
After replacing $g$ with $g'$ as in Remark~\ref{re:decomp},
we can assume $g\in \Or (X,V\oplus \Bc)$.
Whenever $\|\Im f(x)\|\ge 1,$ we have 
\[
\|\Im g(x)\| > 2\|\Im f(x) \|-1 \ge \|\Im f(x)\|.
\]
Whenever $\|\Im f(x) \|< 1,$ we have 
\[
\|\Im g(x)\|> 2\|\Im x \|/(\|\Im a \| -r)) \|-1 
\ge 1>\|\Im f(x) \|
\] 
on $B(a;\mu r).$
\endproof

Now we are ready to prove Proposition~\ref{pr:descent} (and then 
Theorems~\ref{th:gen_nbhdbasis} and \ref{th:nbhdbasis} will be fully proved).
\proof[Proof of Proposition~\ref{pr:descent}]
With $\mu$ as in Hypothesis~\ref{hyp:Runge}, let $0<2^{8}\alpha<\mu^2<1$.
By Proposition~\ref{pr:cover}, each $\PN, N=1,2,\ldots$ has a finite cover
 $\{B_k\}_{k=1}^n$, $B_k=B(x_k,r_k)$ such that 
$B(x_k;2r_k/\mu) \in \mathfrak{B}_P.$
For each $k$, there are a complexified
Banach space $V_k$ and a
\[ 
g_k \in \Or (B(x_k;2r_k/\mu) \cup \c (B(x_k;2r_k/\mu)) ,V_k)
\]
such that
$u \le \|\Im f_k\|$ on $B(x_k;2r_k/\mu) \cup \c (B(x_k;2r_k/\mu))$.
  By Proposition~\ref{pr:balldom}, there are a
complexified Banach space $E_k$ and an $h_k \in \Or (X,E_k)$ such that
$\|\Im h_k \| \ge \|\Im g_k \|$ on $B_k$.  Let
$W_N=\overline{\bigoplus} V_k$, and let $\phi_N \in \Or (X,W_N)$ be defined by
$\phi_N(x)=(h_1(x),h_2(x),\ldots ,h_n(x))$.
Then $\|\Im \phi_N \| \ge u$ on $\PN$.
By Proposition~\ref{pr:tech2}, there is a complexified Banach space $Z_N$ and
$f_N \in \Or (P,Z_N)$ such that
\renewcommand\theenumi{\roman{enumi}}
\renewcommand\labelenumi{(\theenumi)}
\begin{enumerate}
\item $\| f_N(x)\| \le 1$ if $x\in P_N\langle\beta\rangle$, and
\item $\| \Im f_N(x)\|\ge\| \Im h(x)\|$ if
$x\in P_{N+1}\langle\alpha\rangle\setminus P_N\langle\alpha\rangle$.
\end{enumerate}
Now we can take $Z_0 = W_1,$ $f_0 = \phi_1,$
$Y=\overline{\bigoplus}_{N=0}^\infty Z_N,$
and
$f=(f_0,f_1,f_2,\ldots)$.  By Proposition~\ref{pr:ballprops},
for each $x \in P$ there is a neigborhood of $x$ contained in some $\PN$.
Therefore the sequence $\{f_N\}$ is locally bounded, so $f\in \Or (P,Y),$ and
$\|\Im f\| \ge u$.
\endproof

\section{Acyclicity}\label{se:acyclicity}
In order to prove Theorem~\ref{th:acyclic}, we will require the following
technical topological proposition, whose proof is similar to that
of \cite[Proposition 2]{Cartan}.
\begin{prop}\label{pr:top} 
Let $X$ be a paracompact Hausdorff space, $A \subset X$ a closed subset,
$\{U_i\}_{i\in I}$ a relatively open cover of $A$, and for each $i\in I,$
$\tilde{U}_i\subset X$ a neighborhood of $U_i$.   
Let $q$ be a positive integer. 
For each $I' \subset I$ of cardinality at most $q,$ let a  
neighborhood $\tilde{U}_{I'} \subset X$ of 
$\bigcap_{i \in I'} U_i$ be given.  Then there are a neighborhood $P$ of $A,$
a function $\sigma : J \rightarrow I$, and an open cover 
$\{V_j\}_{j \in J}$
such that $\overline{V}_j \subset \tilde{U}_{\sigma j},$ and 
$\bigcap_{j\in J'} \overline{V}_j \cap P \subset \tilde{U}_{\sigma J'}$ 
for all $J' \subset J$ of cardinality at most $q$.
\end{prop}

\begin{proof}
By passing to a refinement, we may assume that $\{\tilde{U}_i\}$ is locally 
finite.  Choose any locally finite open cover $\{V_j\}_{j\in J}$ of $A$
and $\sigma :J\rightarrow I$ so that 
$\overline{V}_j \subset \tilde{U}_{\sigma j}$.
Then for each $x \in A$ there is a neighborhood $P(x)$ of $x$ such that
\renewcommand\theenumi{\roman{enumi}}
\renewcommand\labelenumi{(\theenumi)}
\begin{enumerate}
\item $P(x) \cap V_j \ne \emptyset$ if and only if 
$x \in \overline{V}_j$, and
\item $P(x)$ is contained in the intersection of all $\tilde{U}_{I'}$
 containing $x$
(of which there are finitely many).
\end{enumerate}
Now define $P$ as the union of all $P(x).$  Given $J' \subset J$ of cardinality
at most $q,$ and $y \in \bigcap_{j \in J'} \overline{V}_j \cap P,$ we must show
that $y \in U_{\sigma J'}$.  Since $y \in P,$ $y \in P(x)$ for some $x \in A.$
But $P(x) \cap V_j \ne \emptyset$ for all $j \in J',$ so 
$x \in \bigcap_{j\in J'} \overline{V}_j \cap A \subset \tilde{U}_{\sigma J'}.$
Therefore, $P(x) \subset \tilde{U}_{\sigma J'}$.
\end{proof}
\begin{lem}\label{le:bundle}
Let $\Xr$ be a real Banach space, $\Omega \subset \Xr$ open, and 
$F \rightarrow \Omega$ a real analytic Banach bundle. There exist a 
neighborhood $W \subset X$ of $\Omega$ and a holomorphic Banach bundle
$E \rightarrow W$ whose restriction to $\Omega$ is the complexification
$F \otimes \Bc$ of $F$.
\end{lem}
\begin{proof}
We can assume that $\Omega$ is connected, in which case all fibers of $F$ are
isomorphic.  As in the finite dimensional case, $F$ is determined by 
an open cover $\{U_j\}$ together with a real Banach space $\Yr$ and real 
analytic transition functions
\[
g_{ij}:U_i \cap U_j \rightarrow  \End (\Yr)
\]
satisfying
\renewcommand\theenumi{\roman{enumi}}
\renewcommand\labelenumi{(\theenumi)}
\begin{enumerate}
\item $g_{ij}(x) g_{ji}(x)=id_{\Yr}$ for $x \in U_i \cap U_j$, and
\item $g_{ij}(x) g_{jk}(x) g_{ki}(x)=id_{\Yr}$ for 
$x\in U_i \cap U_j \cap U_k$.
\end{enumerate}
It is not hard to see that  the complexification of $\End (\Yr)$ can be 
naturally
identified with $\End (Y).$  By the definition of real analytic  functions,
each $g_{ij}$ can
be holomorphically extended to a function $\tilde{g}_{ij}$ on a neighborhood 
$\tilde{U}_{ij} \subset X$ of $U_i \cap U_j$.  We can take $\tilde{U}_{ij}$ 
so small that
$\tilde{g}_{ij}(x) \tilde{g}_{ji}(x)=\text{id}_{Y}$ on $\tilde{U}_{ij}$.
Further, choose neighborhoods 
$\tilde{U}_{ijk}$ of $U_i \cap U_j \cap U_k$ so that
$\tilde{g}_{ij}(x) \tilde{g}_{jk}(x) \tilde{g}_{ki}(x)=\text{id}_{Y}$ for 
$x\in \tilde{U}_i \cap \tilde{U}_j \cap \tilde{U}_k$.
If $\Omega$ is an open subset of $\Xr,$
then by Proposition~\ref{pr:top}, after replacing 
$\{U_j\}$ with a refinement, we can find a neighborhood $V_j\subset X$ of 
each $U_j$ such that $V_i \cap V_j \subset \tilde{U}_{ij}$. Then $g_{ij}$
extends to $V_i \cap V_j$ holomorphically.  The holomorphic extensions 
$\tilde{g}_{ij}$ define a holomorphic vector bundle on $W=\bigcup V_i.$
\end{proof}
We will now prove Theorem~\ref{th:acyclic}; in fact, we will prove the
following more general version:
\begin{thm}\label{th:genacyclic}
Let $\Xr$ be a real Banach space with a Schauder basis satisfying
Hypothesis~\ref{hyp:Runge}, let
$\Omega \subset \Xr$ be open, $F \rightarrow \Omega$ a real analytic Banach
bundle, and $\Ca^r$ the sheaf of real analytic $r-$forms with values
in $F$.  Then $H^p (\Omega, \Ca^r)=0$ for $p \ge 1, \, r \ge 0.$
\end{thm}
\begin{proof}
Let $\mfU$ be an open cover of $\Omega$.
Consider a real analytic $p$-cocycle $c\in C^p(\mfU,\Ca^r),$ $p \ge 1$.
We wish to show that, after sufficient refinement of $\mfU$, $c$ becomes a
coboundary.  We accomplish this by complexification.  It can be assumed that
$F$ is trivial over each $U \in \mfU$.  Let $E \rightarrow W$ be
as in Lemma~\ref{le:bundle}.
 
We extend each component of $c$ on
$\Omega$
 to some holomorphic section of $E$ over some neighborhood 
$\tilde{U}_{I'} \subset X$ of each $p+1$-fold intersection 
$\bigcap_{i\in I'} U_i$, and 
construct
 the corresponding neigborhood $P$ of $\Omega$ and open cover 
$\mfV=\{V_j\}$, as 
in Proposition~\ref{pr:top}.  In view of Theorem~\ref{th:nbhdbasis}, we can 
take $P$ and each $V_j$ to be pseudoconvex.  

This enables us to apply the following theorem from \cite{Lempert:vanishing}:
\begin{thm}\label{th:vanishing}
Suppose $X$ is a Banach space with a Schauder basis and 
Hypotesis~\ref{hyp:Runge} holds.  If $P \subset X$ is open and
 pseudoconvex, $E\rightarrow P$ a locally trivial holomorphic
Banach bundle, and $q\ge 1,$ then $H^q(P,E)=0.$
\end{thm}
This implies that  we can find a holomorphic
 cochain $b \in C^{p-1}(\mfV, E)$ whose coboundary is the extension of $c$.
Taking fiberwise the real part of $b|_\Omega$ we obtain a 
$b' \in C^{p-1}(\mfU,F)$ with $\delta b'=c$.
\end{proof}
Theorem~\ref{th:iso} follows, again in greater generality:
\begin{thm}\label{th:geniso}
If $X$ is a complex Banach space with a Schauder basis, and satisfies
Hypothesis~\ref{hyp:Runge}, then
\[
H^q(\Omega,E) \approx
\frac{\text{Ker }\{\dbar : \Gamma (\Omega, \Ca_{0,q})   \rightarrow
\Gamma (\Omega, \Ca_{0,q+1})\}}
     {\text{Im } \{\dbar : \Gamma (\Omega, \Ca_{0,q-1}) \rightarrow
\Gamma (\Omega, \Ca_{0,q})  \}}.
\]
\end{thm}
\begin{proof}
First, recall the local solvability of the $\overline{\partial}-$equation
mentioned in Section~\ref{se:intro} and proved in 
\cite[Proposition 3.2]{Lempert:DolbeaultI}.
Second, observe that $\Ca^r=\bigoplus_{p+q=r} \Ca_{p,q}.$  Then
$H^n(\Omega,\Ca^r)=0$ implies $H^n(\Omega, \Ca_{0,r})=0$
for all $r \ge 0,\, n \ge 1.$ These are the two ingredients required in
the hypothesis of the  abstract 
deRham Theorem (see, for example, \cite[Theorem 3.13]{Wells}).  The isomorphism
theorem follows at once from this.
\end{proof}

\end{document}